\documentclass[11pt]{amsart}
\usepackage[a4paper,margin=1.08in]{geometry}
\usepackage{amsmath,amssymb,amsthm,mathtools}
\usepackage{microtype}
\usepackage[hidelinks]{hyperref}
\hypersetup{pdftitle={Lipschitz chain approximation of integral currents in CAT(0) spaces}}
\allowdisplaybreaks[1]

\newtheorem{theorem}{Theorem}[section]
\newtheorem{lemma}[theorem]{Lemma}
\newtheorem{proposition}[theorem]{Proposition}
\newtheorem{corollary}[theorem]{Corollary}
\theoremstyle{definition}

\newcommand{\M}{\mathbf M}
\newcommand{\I}{\mathbf I}
\newcommand{\N}{\mathbf N}
\DeclareMathOperator{\Vol}{Vol}
\DeclareMathOperator{\spt}{spt}
\DeclareMathOperator{\Lip}{Lip}
\DeclareMathOperator{\car}{car}
\DeclareMathOperator{\Hull}{Hull}
\newcommand{\R}{\mathbb R}
\newcommand{\Z}{\mathbb Z}
\newcommand{\Nbh}{\mathcal N}
\newcommand{\ga}{\gamma}
\newcommand{\llb}{[\![}
\newcommand{\rrb}{]\!]}

\title[Lipschitz chain approximation in CAT(0) spaces]{Lipschitz chain approximation of integral currents in CAT(0) spaces}
\author{Stephan Stadler}
\date{\today}
\begin{document}

\begin{abstract}
We establish a Federer--Fleming type approximation theorem for
integral currents in CAT(0) spaces, without assumptions
on dimension or properness. Every integral current can be
approximated in normal mass by Lipschitz chains whose parametrized
volumes, together with those of their boundaries, converge to
the corresponding current masses. If the boundary
is already represented by a Lipschitz cycle, the approximation can
preserve this cycle exactly. We also lift currents from ultralimits
with prescribed boundaries and converging mass. As a consequence,
flat limits of minimizing currents with uniformly bounded normal
mass are minimizing. 
\end{abstract}
\maketitle

\section{Introduction}
The Federer--Fleming deformation theorem is a fundamental tool in
geometric measure theory. Among its consequences is the approximation
of integral currents by Lipschitz chains in normal mass
\cite[Theorems~5.5 and~5.8]{FF}. It is natural to ask how far this approximation principle extends to general metric spaces.

Basso, Wenger and Young established normal-mass approximation
in complete metric spaces with finite Nagata dimension and suitable
Lipschitz connectivity assumptions \cite[Corollary~1.5]{BWY}.
In another direction,
Goldhirsch obtained normal-mass approximation for compactly
supported integral currents in locally compact metric spaces
satisfying a local Lipschitz extension property
\cite[Theorem~1.3]{Goldhirsch}.
The dimensional restriction in the former result and the local
extension hypothesis in the latter are not available for
arbitrary CAT(0) spaces.

The purpose of this paper is to establish the approximation
theorem for all CAT(0) spaces, without assumptions
on properness or dimension. 

We write $\I_k(X)$ for the group of integral $k$-currents
in $X$, $\M$ for mass, and
$\N(T)=\M(T)+\M(\partial T)$ for normal mass.
For a Lipschitz chain $\alpha$, we denote by $[\alpha]$
the induced integral current and by $\Vol(\alpha)$
its parametrized volume.

\begin{theorem}\label{thm:main}
Let $X$ be a CAT(0) space and let
$T\in\I_k(X)$, $k\geq1$.
There are Lipschitz $k$-chains $\alpha_i$ in $X$
such that
\[
 \N(T-[\alpha_i])\longrightarrow0,\qquad
 \Vol(\alpha_i)\longrightarrow\M(T),\qquad
 \Vol(\partial\alpha_i)\longrightarrow\M(\partial T).
\]
\end{theorem}

In fact, we prove a stronger relative approximation theorem.
It allows us to preserve a prescribed Lipschitz boundary
exactly, as a singular chain and not just as an integral current.

\begin{theorem}\label{thm:strong}
Let $X$ be a  CAT(0) space,
$T\in\I_k(X)$ with $k\geq1$, and let
$\beta$ be a Lipschitz $(k-1)$-cycle with
$[\beta]=\partial T$.
For all $\varepsilon,\delta>0$ there is a Lipschitz
$k$-chain $\alpha$ such that
\[
 \partial\alpha=\beta,\qquad
 \M(T-[\alpha])<\varepsilon,\qquad
 \Vol(\alpha)\leq\M(T)+\varepsilon.
\]
Moreover,
\[
 \car\alpha\subset
 \Nbh_\delta(\spt T\cup\car\beta).
\]
\end{theorem}

An immediate consequence is that integral currents and
Lipschitz chains define the same filling volumes.

\begin{corollary}\label{cor:filling-volumes}
Let $X$ be a CAT(0) space and let
$\beta$ be a Lipschitz $(k-1)$-cycle, $k\geq1$.
Then
\[
 \inf_{\substack{\alpha\text{ Lipschitz }k\text{-chain}\\
                  \partial\alpha=\beta}}
       \Vol(\alpha)
 =
 \inf_{\substack{S\in\I_k(X)\\
                  \partial S=[\beta]}}
       \M(S).
\]
\end{corollary}

\begin{proof}
For every Lipschitz $k$-chain $\alpha$ with
$\partial\alpha=\beta$, the induced current
$[\alpha]$ fills $[\beta]$, and
$\M([\alpha])\leq\Vol(\alpha)$.
This gives one inequality. Conversely, apply
Theorem~\ref{thm:strong} to any integral filling
$S$ of $[\beta]$, with arbitrarily small
volume excess.
\end{proof}

Thus passing from integral currents to Lipschitz chains
does not change filling volume. For boundaries of
parametrized Lipschitz spheres of dimension at least two,
recent work of Vikman and Wenger
\cite[Theorem~1.4]{VW} identifies this infimum
with the infimum of parametrized volumes of Lipschitz
ball fillings. In particular, the sharp Euclidean filling
inequality for integral currents in CAT(0) spaces
\cite[Theorem~1.1]{LytchakFilling} transfers to these
Lipschitz filling problems.

Theorem~\ref{thm:main} also identifies the integral
current complex with the completion, in normal mass,
of the Lipschitz chain complex modulo the kernel of
$\alpha\mapsto[\alpha]$.
Indeed, integral currents are complete in normal mass
\cite[Theorem~8.5]{AK}, while integral zero-currents
already correspond to finite integer combinations of
point currents.

% We expect both theorems to extend to CAT$(1)$ spaces. A promising simple approach is to pass to a Euclidean cone over a CAT$(1)$ space, approximate near the unit sphere around the tip and use a  projection.
% In \cite{Goldhirsch} Goldhirsch approximates compactly supported integral currents by Lipschitz chains with respect to normal norm in locally compact metric
% space with {\em local property L}. There are compact CAT(0) spaces without local property~L, cf.~\cite[Theorem~4]{NaorRabani}.

The approximation theorem provides lifts of integral currents
from ultralimits with prescribed boundaries and convergence of mass.
Let $X_i$ be CAT(0) spaces and let $K_i\subset X_i$ be
nonempty compact subsets. Suppose that there are isometric embeddings
$\varphi_i:K_i\to Z$ into a complete metric space $Z$ such that
$\varphi_i(K_i)$ converge in Hausdorff distance to a compact set
$K\subset Z$. Choose $x_i\in K_i$ and let $(X_\omega,x_\omega)$ be
the pointed ultralimit of $(X_i,x_i)$ with respect to a nonprincipal
ultrafilter $\omega$ on $\mathbb N$. The given embeddings induce an
isometric embedding $\iota_\omega:K\to X_\omega$. 
We write $\mathcal F_Z$ for the integral flat norm in $Z$.

\begin{theorem}\label{thm:ultralimit-lifting}
In the setting above, let $m\geq1$ and let $B_i\in\I_{m-1}(K_i)$ be
cycles of uniformly bounded mass such that
$\lim_\omega\mathcal F_Z((\varphi_i)_\# B_i-B)=0$ for some
$B\in\I_{m-1}(K)$. For every $S\in\I_m(X_\omega)$ with
$\partial S=(\iota_\omega)_\# B$, there are $S_i\in\I_m(X_i)$ with
$$
 \partial S_i=B_i,\qquad
 \lim_\omega\M(S_i)=\M(S).
$$
% Moreover, there is a Banach space $E$ containing isometric copies of
% $Z$, $X_\omega$, and all $X_i$ that agree on the prescribed copies
% of $K$ and $K_i$, in which $\lim_\omega\mathcal F_E(S_i-S)=0$.
\end{theorem}

This applies to sequences of minimizing currents.

\begin{corollary}\label{cor:ultralimit-minimizers}
Let $X_i$, $K_i$, $Z$ and $K$ be as above. Suppose $T_i\in\I_m(X_i)$,
$m\geq1$, minimize mass for their boundaries, have uniformly bounded
normal mass, and satisfy $\spt T_i\subset K_i$. Suppose the pushforwards
$(\varphi_i)_\# T_i$ converge in the integral flat norm in $Z$ to
$T\in\I_m(K)$. Then $(\iota_\omega)_\# T$ minimizes mass for its
boundary in $X_\omega$ for every nonprincipal ultrafilter $\omega$, and
$\lim_{i\to\infty}\M(T_i)=\M(T)$.
\end{corollary}

In particular, when the supports of $T_i$ converge in the
Gromov--Hausdorff sense, one can take $K_i=\spt T_i$ and realize
the Hausdorff convergence in a Banach space $Z$. The images of the
supports then have compact closure in $Z$. The normal-mass bound
gives a weakly convergent subsequence by Ambrosio--Kirchheim
compactness and closure \cite[Theorems~5.2 and~8.5]{AK}. The limit
is supported in $K$, and the convergence is flat by
\cite[Theorem~1.4]{WengerFlat}.

Korevaar and Schoen~\cite[Theorem~3.9]{KorevaarSchoen} proved a
related stability result for harmonic maps. 

%The proofs are given in Section~\ref{sec:ultralimits}.

\medskip
\noindent{\em Idea of proof of Theorem~\ref{thm:main}.}
The main difficulty is to correct the boundary of a Lipschitz
approximation without losing control of its volume.
Using a rectifiable representation, we first approximate an
integral current $T$ by a Lipschitz chain $\alpha$  in mass.
This step may create a huge boundary error $\beta-\partial\alpha$.
However, by construction, this error admits a filling of small mass.

The key observation is then that a minimizing filling has
doubling support away from its boundary, by the
monotonicity formula. This provides the finite-dimensional
structure needed for simplicial deformation.
Using a mapping cylinder construction, we arrange that the
prescribed boundary is a Lipschitz neighborhood retract.
We then approximate the minimizing filling outside a thin
boundary collar by the deformation theorem of
Basso--Wenger--Young. A Lipschitz retraction and a short
homotopy restore the boundary exactly, with arbitrarily
small additional volume cost.

Finally, we obtain normal-mass approximation by first
approximating the boundary and then applying the
prescribed-boundary case to the resulting current.
\medskip

\noindent{\bf Acknowledgements.}
It is my pleasure to thank Alexander Lytchak and Stefan Wenger for helpful comments on an earlier version of this paper.
\medskip

\noindent{\bf Use of AI.}
 The mathematical ideas and overall proof strategy of this paper were developed by the author. ChatGPT was used as an auxiliary tool to explore and refine arguments, identify relevant literature, and assist in drafting and improving the exposition. The author takes full responsibility for the mathematical content and the final text.

\section{Preliminaries}

Throughout the paper, all CAT(0) spaces are assumed
to be complete. No properness or finite-dimensionality
is assumed.

\subsection{Currents and Lipschitz chains}\label{subsec:currents}

We use the notation of \cite[Section~3]{LSU}. The subgroup of compactly
supported integral currents is denoted by $\I_{k,c}(X)$. We write $\|T\|$ and
$\spt T$ for the mass measure and support. Open $r$-balls and $r$-neighborhoods are denoted by $B(x,r)$ and
$\Nbh_r(A)=\{x:d(x,A)<r\}$, respectively. We write $T\llcorner A$ for restriction to a
Borel set $A$. For $T\in\I_j(X)$, $j\geq1$, and a
Lipschitz function $f:X\to\R$, the $r$-slices are denoted
 $\langle T,f,r\rangle$.
They satisfy the coarea bound \cite[Theorems~5.6 and~5.7]{AK}
\[
 \int_a^b\M(\langle T,f,r\rangle)\,dr
 \leq\Lip(f)\,\|T\|(\{a<f<b\}).
\]

A Lipschitz $k$-chain is an element of the free abelian group generated by
Lipschitz maps $\sigma:\Delta^k\to X$ on the ordered standard simplex. Its
boundary is the alternating sum of the restrictions to ordered faces. For
$\alpha=\sum_i m_i\sigma_i$ with distinct $\sigma_i$ and $m_i\ne0$, we define its induced integral current, its carrier and its volume by
\[
 [\alpha]=\sum_i m_i(\sigma_i)_\#\llb\Delta^k\rrb,\qquad
 \car\alpha=\bigcup_i \sigma_i(\Delta^k),\qquad
 \Vol(\alpha)=\sum_i |m_i|\Vol(\sigma_i).
\]
Here $\Vol(\sigma_i)$ is parametrized Hausdorff volume, normalized to agree
with Euclidean volume. In degree zero $\Vol(\alpha)=\sum_i|m_i|$. 

By Lytchak's Rademacher theorem \cite[Theorem~1.6]{LytchakDifferentiation}, cf.~\cite[Proposition~11]{StadlerDiscs} and \cite[Proposition~2.3]{LytchakFilling}, a Lipschitz map $f:A\to X$ from a Borel set $A\subset\R^j$ into a CAT(0) space has at almost every $x\in A$ a unique differential $d_xf:\R^j\to T_{f(x)}X$, linear onto a convex Euclidean subspace. It satisfies
\[
 d_{T_{f(x)}X}\bigl(\log_{f(x)} f(x+v),d_xf(v)\bigr)=o(|v|)
 \qquad (v\to0,\ x+v\in A).
\]
Here $T_yX$ is the tangent cone at $y$ and $\log_y z$ is the tangent vector
of the geodesic from $y$ to $z$ of length $d(y,z)$. The metric
differential is $v\mapsto |d_xf(v)|$. The area and mass formulas give
$\M([\alpha])\leq\Vol(\alpha)$ in CAT(0) spaces \cite[Theorem~7]{Kirchheim},
\cite[Theorem~9.5]{AK}.

By the parametric representation theorem \cite[Lemma~4.1 and Theorem~4.5]{AK}, an integer-rectifiable $j$-current $T$, $j\geq1$, can be written as
\[
 T=\sum_i (f_i)_\#\llb\theta_i\rrb,\qquad
 \M(T)=\sum_i\M\bigl((f_i)_\#\llb\theta_i\rrb\bigr).
\]
Here $f_i:K_i\to X$ are bi-Lipschitz maps from compact sets
$K_i\subset\R^j$ with pairwise disjoint images in $\spt T$, and
$\theta_i\in L^1(\R^j,\Z)$ vanishes off $K_i$. In a CAT(0) space the
mass of the $i$-th term is $\int_{K_i}|\theta_i(x)|\operatorname{Jac}(d_x f_i)\,dx$,
with $\operatorname{Jac}$ the Euclidean Jacobian \cite[Theorem~9.5]{AK}.

Let $H:[0,1]\times X\to Y$ be a Lipschitz homotopy and let $H_t=H(t,\cdot)$.
For an integral current $T$, set
$H(T)=H_\#(\llb 0,1\rrb\times T)$.
The product boundary formula \cite[Proposition~3.4]{BWY} gives
\[
 \partial H(T)+H(\partial T)=(H_1)_\#T-(H_0)_\#T.
\]
For a Lipschitz simplex $\sigma:\Delta^j\to X$, define $H(\sigma)$ by
composing $(t,u)\mapsto H(t,\sigma(u))$ with the affine parametrizations of
the ordered simplices
\[
 [v_0^0,\ldots,v_i^0,v_i^1,\ldots,v_j^1],\qquad 0\leq i\leq j,
\]
with signs $(-1)^i$, where $v_i^t=(t,v_i)$ for the ordered vertices
$v_0,\ldots,v_j$ of $\Delta^j$. Extend this operation linearly to chains.
The prism identity \cite[proof of Theorem~2.10]{Hatcher} gives
\[
 \partial H(\alpha)+H(\partial\alpha)=(H_1)_\#\alpha-(H_0)_\#\alpha.
\]
The signed simplices above induce the product current
$\llb 0,1\rrb\times\llb\Delta^j\rrb$. Pushing forward gives
$[H(\alpha)]=H([\alpha])$.

If the homotopy has time Lipschitz constant bounded by $h$ and spatial
Lipschitz constant by $L$, then, by \cite[Lemma~3.5]{BWY},
\begin{equation}\label{eq:prism}
 \M(H(T))\leq C_j h L^j\M(T),\qquad
 \Vol(H(\alpha))\leq C_j h L^j\Vol(\alpha).
\end{equation}
For the geodesic homotopy between two Lipschitz maps $f_0$ and $f_1$ in a CAT(0) space one can take
\[
 h=\sup_x d(f_0(x),f_1(x)),\qquad L=\max\{\Lip(f_0),\Lip(f_1)\}.
\]

\subsection{Minimizing fillings and monotonicity}\label{subsec:fillings}

For an integral $j$-cycle $T$, $j\geq1$, in a CAT(0) space,
Wenger's Euclidean isoperimetric inequality \cite[Corollary~1.4]{Wenger} gives
\[
 \inf\{\M(S):S\in\I_{j+1}(X),\ \partial S=T\}
 \leq C_j\M(T)^{(j+1)/j}.
\]
For compactly supported $T$, the infimum is attained by a compactly supported minimizing
current \cite[Theorem~1.6]{Wenger}.

Let $P=\sum_{i=1}^N(\llb q_i\rrb-\llb p_i\rrb)$ be an integral
$0$-current of total coefficient zero, with each point listed
according to its multiplicity. Then
\[
  \inf\{\M(R):R\in\I_1(X),\ \partial R=P\}
  =\min_{\sigma\in\mathfrak S_N}\sum_{i=1}^N d(p_i,q_{\sigma(i)}).
\]
A minimizing filling is obtained by summing the oriented geodesic
segments from $p_i$ to $q_{\sigma(i)}$ for a minimizing permutation
$\sigma$ \cite[Theorems~3.2 and~5.3]{BDP}.

For a minimizing $S\in\I_{j,c}(X)$, $j\geq1$, Wenger's monotonicity formula gives
\[
 \omega_j\leq r^{-j}\|S\|(B(x,r))
 \leq t^{-j}\|S\|(B(x,t))
 \qquad\bigl(x\in\spt S,\ 0<r\leq t<d(x,\spt\partial S)\bigr),
\]
where $\omega_j$ is the volume of the Euclidean unit $j$-ball \cite[Corollary~4.4 and the following remark]{WengerRank}, cf.~\cite[Corollary~3.3]{LSU}. In particular,
\begin{equation}\label{eq:filling-radius}
 d(x,\spt\partial S)\leq\bigl(\M(S)/\omega_j\bigr)^{1/j}
 \qquad (x\in\spt S).
\end{equation}

\begin{lemma}\label{lem:polyhedral-kernel}
Let $P$ be a finite piecewise Euclidean polyhedron
and let $\gamma$ be a Lipschitz $q$-cycle in $P$
such that $[\gamma]=0$.
Then there is a Lipschitz $(q+1)$-chain
$\zeta$ in $P$ with
$\partial\zeta=\gamma$.
If $\dim P\leq q$, one can moreover conclude that
$\Vol(\zeta)=0$.
\end{lemma}

\begin{proof}
A finite piecewise Euclidean polyhedron is locally
Lipschitz contractible. By the comparison theorem
for Lipschitz singular homology and integral current
homology \cite[Theorem~1.3]{Mitsuishi}, the natural map
\[
 H_q^{\mathrm{Lip}}(P;\Z)
 \longrightarrow H_q^{\I}(P)
\]
is an isomorphism.
Since $[\gamma]=0$, the homology class represented
by $\gamma$ vanishes. Thus $\gamma$ is the boundary
of a Lipschitz $(q+1)$-chain $\zeta$.

If $\dim P\leq q$, then every Lipschitz map
from a $(q+1)$-simplex into $P$ has vanishing
$(q+1)$-dimensional parametrized volume.
Hence $\Vol(\zeta)=0$.
\end{proof}

\subsection{Relative Federer--Fleming deformation}

We will use the following  simplicial deformation
theorem for integral currents due to Basso--Wenger--Young \cite[Theorem~A.2]{BWY}. It is a variant of the original Federer--Fleming deformation \cite[Theorem~5.5]{FF} \cite[Theorem~1.1]{White}.

For a subset $A$ in a simplicial complex $\Sigma$ we denote by  $\Hull(A)$ its simplicial hull.

\begin{proposition}\label{prop:simplicial}
Let $\Sigma$ be a finite simplicial complex of dimension at most $n$,
with regular Euclidean simplices of edge length $s>0$. Let $T\in\I_j(\Sigma)$,
$j\geq1$, and let $\alpha$ be a Lipschitz $(j-1)$-cycle with
$[\alpha]=\partial T$. There are a Lipschitz $j$-chain $\beta$ and
an integral $(j+1)$-current $S$ such that
\[
 \partial\beta=\alpha,\qquad T-[\beta]=\partial S,\qquad
 \Vol(\beta)+s^{-1}\M(S)\leq C\bigl(\M(T)+s\Vol(\alpha)\bigr).
\]
Here $C$ depends only on $n,j$, and $\car\beta\cup\spt S$ is contained
in $\Hull(\spt T\cup\car\alpha)$.
\end{proposition}

\begin{proof}[Sketch of Proof]
Let $K=\operatorname{Hull}(\operatorname{spt}T\cup\operatorname{car}\alpha)$.
The classical simplicial deformation of $\alpha$ in $K$ gives a
Lipschitz cycle $\alpha_0$ in $K^{(j-1)}$ and a Lipschitz
$j$-chain $\gamma$ in $K$ with
\[
  \partial\gamma=\alpha_0-\alpha,
  \qquad
  \operatorname{Vol}(\gamma)\le Cs\,\operatorname{Vol}(\alpha).
\]
Set $T_0=T+[\gamma]$. Since $\partial T_0=[\alpha_0]$, we obtain a decomposition
$T_0=P+R+\partial S$ with
$\operatorname{spt}R\subset K^{(j-1)}$ as in \cite[Theorem~A.2]{BWY}. Thus $R=0$, and
\[
  \partial P=[\alpha_0],
  \qquad
  \mathbf M(P)+s^{-1}\mathbf M(S)\le C\mathbf M(T_0).
\]
The theorem also gives $\operatorname{spt}P\cup\operatorname{spt}S\subset K$.

Since $[\alpha_0-\partial P]=0$,
Lemma~\ref{lem:polyhedral-kernel}, applied to
$K^{(j-1)}$, gives a Lipschitz $j$-chain
$\eta$ in $K^{(j-1)}$ such that
$\partial\eta=\alpha_0-\partial P$ and
$\Vol(\eta)=0$.
Set $\beta=P+\eta-\gamma$.
Then $\partial\beta=\alpha$,
$T-[\beta]=\partial S$, and
$\car\beta\subset K$.

By subadditivity of parametrized volume, and since $\eta$ is
carried by the $(j-1)$-skeleton, we have
\[
  \operatorname{Vol}(\beta)
  \le \operatorname{Vol}(P)+\operatorname{Vol}(\eta)
       +\operatorname{Vol}(\gamma)
  = \mathbf M(P)+\operatorname{Vol}(\gamma).
\]
Adding $s^{-1}\mathbf M(S)$ and using the deformation estimate gives
\[
  \operatorname{Vol}(\beta)+s^{-1}\mathbf M(S)
  \le C\mathbf M(T_0)+\operatorname{Vol}(\gamma).
\]
Since $T_0=T+[\gamma]$, we have
\[
  \mathbf M(T_0)
  \le \mathbf M(T)+\mathbf M([\gamma])
  \le \mathbf M(T)+\operatorname{Vol}(\gamma).
\]
Thus the preceding bound is at most
$C\mathbf M(T)+(C+1)\operatorname{Vol}(\gamma)$.
Using $\operatorname{Vol}(\gamma)\le C_0s\operatorname{Vol}(\alpha)$
gives the desired estimate
\[
  \operatorname{Vol}(\beta)+s^{-1}\mathbf M(S)
  \le C'\bigl(\mathbf M(T)+s\operatorname{Vol}(\alpha)\bigr).
\]
\end{proof}

\section{Approximation on doubling supports}

\subsection{Mass approximation in rectifiable charts}

In a first step, we aim to approximate an integral current using charts of a rectifiable representation, cf.~\cite[Proposition~1.4]{Goldhirsch}
and \cite[Lemma~7.5]{BWY}.

For each $t>0$ Euclidean space carries a natural grid of side length $t$ induced by the vertex set $t\Z^j\subset\R^j$. Each cube in the grid inherits an orientation from $\R^j$ and we have a natural chain complex of grid cubes with integer coefficients. Choosing a triangulation, each cubical chain induces a simplicial chain.

\begin{lemma}\label{lem:euclidean-grid}
For $j\geq1$, let 
$\theta\in L^1(\R^j,\Z)$ with compact support $K\subset\R^j$.
For every $t>0$ there is a cubical $j$-chain  $P_t$ of side length $t$.
Its carrier is contained in the union of the closed grid cubes of side length $t$ meeting
$K$, and $[P_t]=\llb\theta_t\rrb$, where
\[
 \|\theta_t-\theta\|_1\longrightarrow0,\qquad
 t\,\Vol(\partial P_t)\longrightarrow0
 \quad (t\to 0).
\]
\end{lemma}

\begin{proof}
To define $P_t$ we choose the multiplicity $m_Q$ of $P_t$ on a grid  cube $Q$ to be the integer nearest to the average $t^{-j}\int_Q\theta$. The cubical chain $P_t$ induces the current $\llb \theta_t\rrb$ such that $\llb \theta_t\rrb\llcorner Q=m_Q\llb Q\rrb$. Its support is a finite union of grid cubes intersecting $K$.
The cube averages converge in $L^1$. Since $\theta$
is integer-valued, rounding increases the $L^1$-error by at most a factor of two and therefore $\|\theta_t-\theta\|_1\longrightarrow0$. 

To estimate the boundary volume, fix a coordinate
direction $e_\ell$. The coefficient of
$\partial P_t$ on the common face of adjacent
cubes $Q$ and $Q+te_\ell$ is, up to orientation,
$m_{Q+te_\ell}-m_Q$.
Summing the contributions of all faces $F$ orthogonal
to $e_\ell$, gives
\[
 t\sum_{F\perp e_\ell}
 \Vol\bigl((\partial P_t)\llcorner F\bigr)
 =
 \|\theta_t(\,\cdot+te_\ell)-\theta_t\|_{L^1}.
\]
Summing over the coordinate directions, we obtain
\[
 t\Vol(\partial P_t)
 =
 \sum_{\ell=1}^j
 \|\theta_t(\,\cdot+te_\ell)-\theta_t\|_{L^1}.
\]
Each summand is bounded by $2\|\theta_t-\theta\|_1+
\|\theta(\,\cdot+t e_\ell)-\theta\|_1$, which tends to zero by continuity of
translations in $L^1$.
\end{proof}

\begin{lemma}\label{lem:charts}
Let $X$ be a CAT(0) space and let $T\in\I_j(X)$, $j\geq1$. For every $\eta>0$ and all sufficiently small $s>0$, there is a Lipschitz $j$-chain $\alpha_s$ such that
\[
 \car\alpha_s\subset\Nbh_s(\spt T),\qquad
 \M(T-[\alpha_s])<\eta,\qquad \Vol(\alpha_s)\leq \M(T)+\eta.
\]
For fixed $\eta$ one can arrange $s\Vol(\partial\alpha_s)\to0$ as $s\to 0$.
\end{lemma}

\begin{proof}
Choose a rectifiable representation
\[
 T=\sum_i(f_i)_\#\llb\theta_i\rrb
\]
as in Section~\ref{subsec:currents}. Select finitely many
charts so that the total mass of the discarded terms
is less than $\eta/2$.
Extend each selected chart $f_i:K_i\to X$
to a Lipschitz map
$\varphi_i:\R^j\to X$ by Kirszbraun's theorem
\cite{LangSchroeder}.

For $s>0$, put
\[
 t_i=\frac{s}{2\sqrt j\max\{1,\Lip(\varphi_i)\}}
\]
and apply Lemma~\ref{lem:euclidean-grid} to
$\theta_i$, obtaining cubical chains $P_{t_i}$.
Set
\[
 \alpha_s=\sum_i(\varphi_i)_\#P_{t_i}.
\]
Every cube occurring in $P_{t_i}$ meets $K_i$
and has diameter $\sqrt j\,t_i$.
Thus $\car\alpha_s\subset\Nbh_s(\spt T)$.
Write $T_i=(f_i)_\#\llb\theta_i\rrb$.
Since $\varphi_i=f_i$ on $K_i$, we have
\[
\begin{aligned}
 \M\bigl(T_i-[(\varphi_i)_\#P_{t_i}]\bigr)
 &\leq
 \Lip(\varphi_i)^j
 \|\theta_i-\theta_{i,t_i}\|_{L^1}
 \longrightarrow0 .
\end{aligned}
\]
Together with the estimate for the discarded charts,
this yields
\[
 \M(T-[\alpha_s])<\eta
\]
for all sufficiently small $s$.

We next estimate the parametrized volume.
By the area formula
\cite[Theorem~7]{Kirchheim}, the volume of the chain
associated with the $i$-th chart satisfies
\[
\begin{aligned}
 \Vol\bigl((\varphi_i)_\#P_{t_i}\bigr)
 &\leq
 \int_{\R^j}
 |\theta_{i,t_i}(x)|
 \operatorname{Jac}(d_x\varphi_i)\,dx\\
 &\leq
 \int_{K_i}
 |\theta_i(x)|\operatorname{Jac}(d_xf_i)\,dx+
 \Lip(\varphi_i)^j
 \|\theta_{i,t_i}-\theta_i\|_{L^1}\\
 &=\M(T_i)+
 \Lip(\varphi_i)^j
 \|\theta_{i,t_i}-\theta_i\|_{L^1}.
\end{aligned}
\]
Here we have used the mass representation formula
\cite[Theorem~9.5]{AK}.
Summing over the finitely many selected charts,
we obtain
\[
 \Vol(\alpha_s)
 \leq \sum_i\M(T_i)+o(1)
 \leq \M(T)+o(1).
\]
Consequently $\Vol(\alpha_s)\leq\M(T)+\eta$
for all sufficiently small $s$.

Finally, the Lipschitz pushforward estimate gives
\[
 \Vol(\partial\alpha_s)
 \leq\sum_i
 \Lip(\varphi_i)^{j-1}\Vol(\partial P_{t_i}).
\]
For each selected chart, $s/t_i$ is independent
of $s$. Hence Lemma~\ref{lem:euclidean-grid}
implies
\[
 s\Vol(\partial\alpha_s)\longrightarrow0.
\]
\end{proof}

\subsection{Boundary correction near doubling sets}

\begin{lemma}\label{lem:mesh-correction}
Let $K$ be a compact doubling subset of a  CAT(0) space $X$,
and let $s>0$. Suppose $T\in\I_{j,c}(X)$, $j\geq1$, and $\alpha$ is a
Lipschitz $(j-1)$-cycle with $[\alpha]=\partial T$. If
$\spt T\cup\car\alpha\subset\Nbh_s(K)$, there are a Lipschitz $j$-chain
$\beta$ and $S\in\I_{j+1,c}(X)$ with
\[
 \partial\beta=\alpha,\qquad T-[\beta]=\partial S,\qquad
 \Vol(\beta)+s^{-1}\M(S)\leq C\bigl(\M(T)+s\Vol(\alpha)\bigr).
\]
Here $\car\beta\cup\spt S\subset\Nbh_{Cs}(\spt T\cup\car\alpha)$, and $C$
depends only on $j$ and the doubling constant of $K$.
\end{lemma}

\begin{proof}
The closed $s$-neighborhood $Y$ of $K$ has a finite cover by closed
$2s$-balls about a maximal $s$-separated set in $K$. Its
$s$-multiplicity is bounded in terms of doubling, since the centers of
balls meeting a set of diameter at most $s$ lie within $5s$ of each other.
By \cite[Lemma~2.1]{LSU} there are a finite simplicial complex $\Sigma$ of bounded dimension metrized as a
subcomplex of some simplex of edge length $s$ in a Euclidean space. Moreover, there are maps $\psi:Y\to\Sigma$ and
$\varphi:\Sigma\to X$ such that $\psi$ is
$C$-Lipschitz, $\varphi$ is $C$-Lipschitz on each simplex, and
$f=\varphi\circ\psi$ is $C$-Lipschitz with $d(y,f(y))\leq Cs$.

Pass to the piecewise Euclidean length metric on $\Sigma$. Then $\psi$ and
$\varphi$ are both locally $C$-Lipschitz. In particular,
\[
 \M(\psi_\#T)\leq C^{j}\M(T),\qquad
 \Vol(\psi_\#\alpha)\leq C^{j-1}\Vol(\alpha).
\]
Apply Proposition~\ref{prop:simplicial} to $\psi_\#T,\psi_\#\alpha$, obtaining
$\tilde\beta$ and $\tilde S$. Let $H$ be the geodesic homotopy from the inclusion of
$Y$ to $f$. Set $\beta=\varphi_\#\tilde \beta-H(\alpha)$ and $S=\varphi_\#\tilde S-H(T)$.
Since $\partial\alpha=0$, the homotopy
identities give
\[
  \partial H(\alpha)=f_\#\alpha-\alpha,\qquad
  \partial H(T)=f_\#T-T-[H(\alpha)].
\]
Consequently, $\partial\beta=\alpha$ and
$\partial S=T-[\beta]$.

The homotopy has time Lipschitz bound $Cs$ and spatial Lipschitz
bound $C$. By \eqref{eq:prism},
\[
  \Vol(H(\alpha))\leq C_{j-1}sC^j\Vol(\alpha),\qquad
  \M(H(T))\leq C_jsC^{j+1}\M(T).
\]
Together with the Lipschitz bound for $\varphi$, this gives
\[
\begin{aligned}
  \Vol(\beta)
  &\leq C^j\Vol(\tilde\beta)+C_{j-1}sC^j\Vol(\alpha),\\
  s^{-1}\M(S)
  &\leq C^{j+1}s^{-1}\M(\tilde S)+C_jC^{j+1}\M(T).
\end{aligned}
\]
We infer
\[
  \Vol(\beta)+s^{-1}\M(S)
  \leq \tilde C\bigl(\M(T)+s\Vol(\alpha)\bigr).
\]

Put $A=\spt T\cup\car\alpha$. By Proposition~\ref{prop:simplicial},
$\car\tilde\beta\cup\spt \tilde S$ lies in $\Hull(\psi(A))$.
Each point $z$ in this hull lies in a simplex containing $\psi(a)$
for some $a\in A$. Hence
\[
  d(\varphi(z),a)
  \leq d(\varphi(z),\varphi(\psi(a)))+d(f(a),a)
  \leq 2Cs.
\]
The homotopy terms are supported in $H([0,1]\times A)$, which lies
within distance $Cs$ of $A$ as claimed.
\end{proof}

\begin{lemma}\label{lem:return}
Let $X$ be a CAT(0) space and $T\in\I_{j,c}(X)$, $j\geq1$, have doubling support. For all $\varepsilon,\delta>0$ there are a Lipschitz
$j$-chain $\alpha$ and a current
$W\in\I_{j,c}(X)$ such that, writing $\beta=\partial \alpha$,
\[
 \Vol(\alpha)\leq\M(T)+\varepsilon,\qquad
 \Vol(\beta)\leq\M(\partial T)+\varepsilon,\qquad
 \partial W=[\beta]-\partial T,\qquad \M(W)<\varepsilon .
\]
One can arrange $\M(T-[\alpha])<\varepsilon$,
$\car\alpha\subset\Nbh_\delta(\spt T)$, and
$\car\beta\cup\spt W\subset\Nbh_\delta(\spt\partial T)$. If
$\partial T=0$, one can take $\beta=W=0$.
\end{lemma}

\begin{proof}
First let $Z$ be a cycle with compact doubling
support. Take a chart approximation $\gamma_s$ from Lemma~\ref{lem:charts}.
Apply Lemma~\ref{lem:mesh-correction} to $Z-[\gamma_s]$ and
$-\partial\gamma_s$ with $K=\spt Z$ to obtain $\lambda$ and $S$.
Then $\beta=\gamma_s+\lambda$ is a cycle and $W=-S$ satisfies
$\partial W=[\beta]-Z$. The chart error and then $s$ can be chosen so that
$\Vol(\beta)-\M(Z)$ and $\M(W)$ are arbitrarily small  and
$\car\beta\cup\spt W$ lies in any prescribed neighborhood of $\spt Z$.
For $Z=\partial T$, the support is a closed subset of the doubling support
of $T$. 

Fix a small $\eta>0$. Choose $\alpha_s$ by Lemma~\ref{lem:charts}, with
\[
 \M(T-[\alpha_s])<\eta,\qquad \Vol(\alpha_s)\leq \M(T)+\eta.
\]
Take $s$ small enough such
that $s\Vol(\partial\alpha_s)<\eta$ and $s(\M(\partial T)+1)<\eta$.
Choose $\beta,W$ by the first step within the $s$-neighborhood of
$\spt\partial T$, with volume excess and $\M(W)$ less than
$\min\{1,\eta\}$. Apply Lemma~\ref{lem:mesh-correction} with $K=\spt T$
to $T+W-[\alpha_s]$ and $\beta-\partial\alpha_s$. The correcting chain
has volume at most $4C\eta$ and its carrier is within $(C+1)s$ of
$\spt T$. Adding it to $\alpha_s$ gives the required chain for small
enough $\eta$ and $s$.
\end{proof}

\section{Approximation with prescribed boundary}

\begin{proposition}\label{prop:filling}
Let $X$ be a CAT(0) space, $T\in\I_k(X)$, $k\geq1$, and
$\beta$ a Lipschitz $(k-1)$-cycle with $[\beta]=\partial T$. For all
$\varepsilon,\delta>0$ there is a Lipschitz $k$-chain $\alpha$ with
$\partial\alpha=\beta$ and $\Vol(\alpha)\leq\M(T)+\varepsilon$. If
$\beta\ne0$, one can arrange
\[
 \car\alpha\subset\Nbh_{C_k\M(T)^{1/k}+\delta}(\car\beta),
\]
where $C_k$ depends only on $k$.
\end{proposition}

\begin{proof}
If $\beta=0$, take $\alpha=0$.
For $k=1$, write
$\beta=\sum_{i=1}^N([q_i]-[p_i])$, listing each point according to
its multiplicity. Choose the indexing so that the matching
$p_i\mapsto q_i$ minimizes $\sum_i d(p_i,q_i)$.
By the least-cost matching formula in Section~\ref{subsec:fillings}, this
sum is the minimum mass of an integral filling of $[\beta]$, hence
is at most $\M(T)$. Let $\alpha$ be the sum of the geodesic
segments oriented from $p_i$ to $q_i$. Then
$\partial\alpha=\beta$ and
\[
  \Vol(\alpha)\leq\sum_{i=1}^N d(p_i,q_i)\leq\M(T).
\]
Each segment starts in $\car\beta$ and has length at most $\M(T)$,
so $\car\alpha\subset\Nbh_{\M(T)+\delta}(\car\beta)$.

\emph{Lifting the boundary.}
Suppose $k\geq2$. Choose a simplicial complex $\Sigma$, realized as a subcomplex of a Euclidean
simplex in $\R^N$,
and a Lipschitz map
$f:\Sigma\to X$ with $f_\#\ga=\beta$ where $\ga$ denotes a fundamental class, cf.~\cite[Section~4.1]{VW}. 

Let $Y=X\times\R^N$, let $\iota:X\to X\times\{0\}$ be the  inclusion, and let $\pi:Y\to X$ be projection. For $0<s\leq1$ set
\[
 g_s(u)=(f(u),su),\qquad \Gamma=g_s(\Sigma),\qquad
 \beta_s=(g_s)_\#\gamma .
\]
The graph $\Gamma$ is bi-Lipschitz to $\Sigma$, and $\pi_\#\beta_s=\beta$. The
homotopy $H(t,u)=(f(u),tsu)$ defines a Lipschitz chain $H_s$ 
with $\partial H_s=\beta_s-\iota_\#\beta$ and $\Vol(H_s)\leq C_\beta s$
by \eqref{eq:prism}. The current $\iota_\#T+[H_s]$ fills $[\beta_s]$.
Thus there is a compactly supported minimizing filling $S$ in $Y$ with
$\M(S)\leq\M(T)+C_\beta s$.

\emph{Cutting off a collar.}
We use a Lipschitz retraction $\rho$ onto $\Gamma$ near the graph. Let
$r_s$ be a Lipschitz neighborhood retraction onto $s\Sigma$ in $\R^N$
\cite[Section~3]{Goldhirsch}. Then $\rho(x,v)=g_s(r_s(v)/s)$ retracts
$\Nbh_{r_0}(\Gamma)$ onto $\Gamma$ for some $r_0>0$, with
\[
 d(y,\rho(y))\leq(1+\Lip(\rho))d(y,\Gamma).
\]

Put $d_\Gamma=d(\,\cdot\,,\Gamma)$. For almost every $r>0$ the currents
$U_r=S\llcorner\{d_\Gamma<r\}$ and $V_r=S-U_r$ are integral. Set
$A_r=\partial V_r=-\langle S,d_\Gamma,r\rangle$. Then
$\spt A_r\subset\{d_\Gamma=r\}$ and $\partial U_r=[\beta_s]-A_r$.
The support of $V_r$ lies in $K_r=\spt S\cap\{d_\Gamma\geq r\}$. 
Since $S$ is minimizing, monotonicity implies for every $x\in K_r$ and $0<u\leq r/2$,
\[
  \omega_k u^k \leq \|S\|(B(x,u)) \leq C u^k,
  \qquad C=(2/r)^k\M(S).
\]
Thus the compact set $K_r$ is doubling.
The coarea inequality gives
\[
 \int_t^{2t}\M(A_r)\,dr\leq\|S\|(\{t<d_\Gamma<2t\})\longrightarrow0
 \qquad(t\to 0).
\]
Consequently there are good radii $r\to 0$ with $r\M(A_r)\to0$.

\emph{Restoring the boundary.}
Restrict to such good radii $r<r_0/2$. Apply Lemma~\ref{lem:return} to $V_r$,
with error $\eta>0$ and $\delta$ less than $r$.
This gives a Lipschitz $k$-chain $\alpha_r$ with
$\gamma_r=\partial\alpha_r$ and an integral $k$-current $W$ such that
\[
  \Vol(\alpha_r)\leq\M(V_r)+\eta,\qquad
  \Vol(\gamma_r)\leq\M(A_r)+\eta,\qquad
  \partial W=[\gamma_r]-A_r.
\]
Moreover, $\car\alpha_r\subset\Nbh_r(\spt V_r)$ and
$\car\gamma_r\cup\spt W\subset\Nbh_r(\spt A_r)$.
Since $\spt A_r\subset\{d_\Gamma=r\}$, the latter neighborhood lies
in $\Nbh_{2r}(\Gamma)$.

Both $U_r$ and $W$ are supported in the domain of $\rho$, and
$\rho_\#U_r=\rho_\#W=0$ since $\dim\Gamma\leq k-1$.
Pushing forward $\partial U_r=[\beta_s]-A_r$ and
$\partial W=[\gamma_r]-A_r$, and using that $\rho$ fixes $\Gamma$,
gives
\[
  [\beta_s]=\rho_\#A_r=[\rho_\#\gamma_r].
\]
Hence,
$\beta_s-\rho_\#\gamma_r$ bounds a Lipschitz $k$-chain $\zeta$
in $\Gamma$. Since $\dim\Gamma\leq k-1$, we have $\Vol(\zeta)=0$.

Let $H$ be the geodesic homotopy from the inclusion of
$\Nbh_{2r}(\Gamma)$ to $\rho$. Since $\partial\gamma_r=0$,
\[
  \partial H(\gamma_r)=\rho_\#\gamma_r-\gamma_r.
\]
Thus $\widetilde\alpha=\alpha_r+H(\gamma_r)+\zeta$ satisfies
$\partial\widetilde\alpha=\beta_s$.
For $y\in\car\gamma_r$, we have
$d(y,\rho(y))\leq 2(1+\Lip(\rho))r$. The spatial Lipschitz bound
of $H$ is at most $\max\{1,\Lip(\rho)\}$, so \eqref{eq:prism} gives
\[
  \Vol(H(\gamma_r))
  \leq C_\rho r\Vol(\gamma_r)
  \leq C_\rho r\bigl(\M(A_r)+\eta\bigr).
\]
Consequently,
\[
  \Vol(\widetilde\alpha)
  \leq \M(V_r)+\eta+C_\rho r\bigl(\M(A_r)+\eta\bigr)
  \leq \M(T)+C_\beta s+\eta+C_\rho r\bigl(\M(A_r)+\eta\bigr).
\]

The filling radius bound \eqref{eq:filling-radius} implies that $\spt V_r$
lies at  distance at most $(\M(S)/\omega_k)^{1/k}$ from $\Gamma$ and therefore
$\car\alpha_r$ lies at distance at most
$(\M(S)/\omega_k)^{1/k}+r$ from $\Gamma$. Each point of
$\car H(\gamma_r)$ lies on a geodesic from some
$y\in\car\gamma_r$ to $\rho(y)\in\Gamma$. Its distance to
$\Gamma$ is bounded by its length, which is at most
$2(1+\Lip(\rho))r$. Since $\car\zeta\subset\Gamma$, the entire
carrier of $\widetilde\alpha$ lies within distance
\[
  (\M(S)/\omega_k)^{1/k}+2(1+\Lip(\rho))r
\]
of $\Gamma$.

We now choose the parameters. First choose $s>0$
so small that the error $C_\beta s$
in the filling mass is arbitrarily small.
For this fixed $s$, the neighborhood retraction
$\rho$ is fixed as well. In particular,
the constants $\Lip(\rho)$ and $C_\rho$
may depend on $s$, but not on $r$.

Next choose $r$ along the good radii tending
to zero, so that
\[
 r\M(A_r)\longrightarrow0,\qquad
 r(1+\Lip(\rho))\longrightarrow0,\qquad
 C_\rho r\M(A_r)\longrightarrow0.
\]
Finally, choose $\eta>0$ sufficiently small.
The preceding estimates then give
\[
 \Vol(\widetilde\alpha)\leq\M(T)+\varepsilon,
\]
and, after making the choices sufficiently small,
\[
 \car\widetilde\alpha
 \subset
 \Nbh_{C_k\M(T)^{1/k}+\delta}(\Gamma).
\]
The projection $\pi:Y\to X$ is $1$-Lipschitz,
sends $\Gamma$ onto $\car\beta$, and satisfies
$\pi_\#\beta_s=\beta$.
Consequently
 $\alpha=\pi_\#\widetilde\alpha$
has all required properties.
\end{proof}

\begin{corollary}\label{cor:optimal-radius}
Let $X$ be a CAT(0) space and let
$\beta\ne0$ be a Lipschitz $(k-1)$-cycle,
$k\geq1$. Set
\[
\operatorname{Fill}(\beta)=\inf\{\M(S):S\in\I_k(X),\ \partial S=[\beta]\}.
\]
Then for every $\varepsilon,\delta>0$
there is a Lipschitz $k$-chain $\alpha$
with
\[
 \partial\alpha=\beta,\qquad
 \Vol(\alpha)\leq \operatorname{Fill}(\beta)+\varepsilon,
\]
and $\car\alpha\subset
 \Nbh_{C_k \operatorname{Fill}(\beta)^{1/k}+\delta}(\car\beta)$.

\end{corollary}

% \begin{proof}
% Choose an integral filling $T$ of $[\beta]$
% with $\M(T)$ arbitrarily close to $F$
% and apply Proposition~\ref{prop:filling}.
% The assertion follows by continuity of
% $t\mapsto t^{1/k}$.
% \end{proof}

\begin{proof}[Proof of Theorem~\ref{thm:strong}]
Use Lemma~\ref{lem:charts} to choose $\alpha_0$ carried in
$\Nbh_{\delta/3}(\spt T)$, with volume at most $\M(T)+\tau$ and
$\M(T-[\alpha_0])<\tau$. Apply Proposition~\ref{prop:filling} to
$T-[\alpha_0]$ and $\beta-\partial\alpha_0$. It gives a chain $\gamma$
with boundary $\beta-\partial\alpha_0$, volume less than $2\tau$, and
carrier within $C_k\tau^{1/k}+\delta/3$ of the carrier of that cycle.
If the cycle is zero, take $\gamma=0$. Then $\alpha=\alpha_0+\gamma$ has
boundary $\beta$, mass error less than $3\tau$, and volume at most
$\M(T)+3\tau$. Take $3\tau<\varepsilon$ and $C_k\tau^{1/k}<\delta/3$.
\end{proof}

\begin{proof}[Proof of Theorem~\ref{thm:main}]
For $k=1$, apply Theorem~\ref{thm:strong} with the point chain inducing
$\partial T$. For $k\geq2$, apply it first to $\partial T$ with
prescribed boundary zero. This gives Lipschitz cycles $\beta_i$ with
$\M(\partial T-[\beta_i])\to0$ and $\Vol(\beta_i)\to\M(\partial T)$.
By the Euclidean isoperimetric inequality, $\partial T-[\beta_i]$ bounds
a current $W_i$ with $\M(W_i)\to0$. Apply Theorem~\ref{thm:strong} to
$T-W_i$ with boundary $\beta_i$ and error $1/i$. Then
\[
 \M(T-[\alpha_i])\leq \M(W_i)+1/i,\qquad
 \M([\alpha_i])\leq\Vol(\alpha_i)\leq \M(T)+\M(W_i)+1/i .
\]
This proves the mass and volume convergence, and $\partial\alpha_i=\beta_i$ gives the boundary convergence.
\end{proof}

\section{Lifting currents from ultralimits}\label{sec:ultralimits}

Wenger's estimate for the lifting of disc fillings
\cite[Theorem~5.1]{WengerAlmostEuclidean} is a two-dimensional
precedent. Here we lift Lipschitz simplices with prescribed faces
and correct the boundary of the resulting currents. To compare
the lifted currents with their limits, we use construction of Ikonen and Wenger
\cite[Theorem~1.6 and Section~4]{IkonenWenger}.

\begin{lemma}\label{lem:lift-simplex}
Let $X_i$ and $X_\omega$ be CAT(0) spaces isometrically
embedded in a Banach space $E$. Let $f:\Delta^d\to X_\omega$ be
Lipschitz and suppose that
$\lim_\omega d_E(z,X_i)=0$ for every $z\in f(\Delta^d)$. Let
$g_i:\partial\Delta^d\to X_i$ be uniformly Lipschitz maps converging
uniformly in $E$ along $\omega$ to $f|_{\partial\Delta^d}$.
There are uniformly Lipschitz extensions $f_i:\Delta^d\to X_i$ of
$g_i$ converging uniformly to $f$ with respect to $\omega$ and satisfying
$\limsup_\omega\Vol(f_i)\leq\Vol(f)$.
\end{lemma}

\begin{proof}
For $d=0$, lift the image point. For $d\geq1$, choose increasingly
fine finite nets of interior points of $\Delta^d$, and lift their
images. At each fixed finite net, $g_i$  extends to the net with a uniform Lipschitz constant.
The Lang--Schroeder extension theorem \cite{LangSchroeder} and a
diagonal argument give uniformly Lipschitz extensions $u_i$ of $g_i$ converging
uniformly to $f$.

Fix $0<\varepsilon<1$. Using metric differentials
as in
\cite[Proposition~4.5]{WengerAlmostEuclidean}, we obtain finitely many
pairwise disjoint compact sets $K_a\subset\operatorname{int}\Delta^d$
covering all but $\varepsilon$ of the Lebesgue measure of $\Delta^d$,
and Euclidean norms $\|\cdot\|_a$ such that
\[
 d(f(x),f(y))\leq \|x-y\|_a\quad(x,y\in K_a),\qquad
 \operatorname{Jac}(\|\cdot\|_a)\leq \operatorname{Jac}(|d_xf|)+\varepsilon
 \quad\text{for a.e. }x\in K_a .
\]
% Indeed, apply the decomposition to the bi-Lipschitz graph
% $x\mapsto(f(x),b x)$ in the Euclidean product, with $b>0$ small.
% Its metric differentials are Euclidean norms bounded below by
% $b|\cdot|$. Subdivide the pieces by diameter and by the values of the
% differential to obtain distance comparison with a fixed Euclidean
% norm on each piece, with arbitrarily small multiplicative error.
% Absorb the upper error into $q_a$ and retain finitely many compact
% pieces by inner regularity. The graph Jacobians converge uniformly
% to the original Jacobians as $b\downarrow0$, since they are products
% of $\sqrt{s_\ell^2+b^2}$ with the singular values $s_\ell$ bounded
% by $\Lip(f)$. The norms $q_a$ have a common upper bound in the
% Euclidean metric, independent of $\varepsilon$.

Now lift $f|_{K_a}$ using finite nets and the extension theorem for the
Euclidean metric given by $\|\cdot\|_a$. This gives maps $v_{a,i}:K_a\to X_i$
converging uniformly to $f|_{K_a}$ with respect to $\omega$, with $\Lip(v_{a,i})\leq 1+o(1)$. 
Choose $0<h<\varepsilon$ so that the $h$-neighborhoods of the $K_a$
are pairwise disjoint and avoid $\partial\Delta^d$. Let $U$ be the
complement of these neighborhoods in $\Delta^d$. On
$U\cup\bigcup_a K_a$, prescribe a map $w_i$ by setting $w_i=u_i$ on
$U$ and $w_i=v_{a,i}$ on $K_a$. Since both $u_i$ and $v_{a,i}$
converge uniformly to $f$ on $K_a$ along $\omega$, we have
$d(w_i,u_i)\leq h$ on the prescribed domain for $\omega$-almost all $i$.

On $U$, the Lipschitz bound for $w_i$ is that of $u_i$. On each
$K_a$, the bound $\Lip_{(\|\cdot\|_a)}(v_{a,i})\leq 1+o(1)$
gives a uniform Euclidean Lipschitz bound. If $x$ and $y$ lie in different pieces of the
prescribed domain, then $|x-y|\geq h$, and hence
\[
  d(w_i(x),w_i(y))
  \leq d(u_i(x),u_i(y))+2h
  \leq \bigl(\Lip(u_i)+2\bigr)|x-y|.
\]
Thus $w_i$ has a Lipschitz bound independent of $h$ and
$\varepsilon$. Extend $w_i$ to $f_i$ on $\Delta^d$ using
\cite{LangSchroeder}. 
Since $\partial\Delta^d\subset U$, each $f_i$ extends $g_i$.
Every point outside the prescribed domain is within $h$ of some
$K_a$. Comparison at a nearest point of $K_a$ gives
$\sup_{\Delta^d}d(f_i,u_i)\leq Ch$.

On $K_a$, the area formula and the $\|\cdot\|_a$-Lipschitz bound
give
\[
  \limsup_\omega \Vol(f_i|_{K_a})
  \leq \mathcal L^d(K_a)\operatorname{Jac}(\|\cdot\|_a)
  \leq \int_{K_a}
       \bigl(\operatorname{Jac}(|d_xf|)+\varepsilon\bigr)\,dx.
\]
The complement of the $K_a$ contributes at most $C\varepsilon$
by the common Lipschitz bound. Summing gives
$\limsup_\omega\Vol(f_i)\leq\Vol(f)+C'\varepsilon$.
Since $h<\varepsilon$, a diagonal argument as $\varepsilon\to 0$
gives uniform convergence to $f$ with respect to $\omega$ and the desired
volume estimate.
\end{proof}

For an integral cycle $T$ in $Y$, let $\operatorname{Fillvol}_Y(T)$
denote the infimum of the masses of its integral fillings in $Y$.

\begin{lemma}\label{lem:intrinsic-small-fillings}
For each $d\geq0$ there is a constant $C_d$ with the following
property. Let $X$ be a CAT(0) space isometrically
embedded in a Banach space $E$, and let $T\in\I_d(X)$ be a cycle. Then
\[
  \operatorname{Fillvol}_X(T)
  \leq C_d\M(T)\operatorname{Fillvol}_E(T)^{1/(d+1)}.
\]
% In particular, if $X_i$ are CAT(0) spaces isometrically
% embedded in $E$ and $T_i\in\I_d(X_i)$ are cycles of uniformly
% bounded mass with $\operatorname{Fillvol}_E(T_i)\to0$, then
% $\operatorname{Fillvol}_{X_i}(T_i)\to0$.
\end{lemma}

\begin{proof}
We may assume that $E$ is injective, since an isometric embedding
into an injective Banach space can only decrease the ambient
filling volume. Glue two copies $E_0,E_1$ of $E$ to
$X\times[0,1]$ along $X\times\{0\}$ and $X\times\{1\}$,
respectively.

The resulting space $Y$ is complete with respect to the induced length metric.
It is not hard to show that it supports cone inequalities in all degrees since the radial contractions in $X$ extend to contractions on $Y$ by injectivity of $E$. Thus $Y$ satisfies Euclidean isoperimetric inequalities in all degrees by \cite[Theorem~1.2]{Wenger}.

Now let $S\in\I_{d+1}(E)$ satisfy $\partial S=T$, and put
$s=\M(S)$. Let $S^j$ and $T^j$ be the copies of $S$ and $T$ in
$E_j$. In $Y$ the current
\[
  U=S^0+\llb 0,1\rrb\times T-S^1
\]
is a cycle since the cylinder has boundary $T^1-T^0$. Its
mass is at most $2s+c_d\M(T)$. Thus there is a filling
$W\in\I_{d+2}(Y)$ of $U$ with
\[
  \M(W)\leq C_d\bigl(s+\M(T)\bigr)^{(d+2)/(d+1)}.
\]
Write $\vartheta:Y\to[0,1]$
for the $1$-Lipschitz height function, which is constant on
$E_0$ and $E_1$.
For almost every $t\in(0,1)$, the restriction
$U\llcorner\{\vartheta<t\}$ has boundary $T^t$, where $T^t$
denotes the copy of $T$ on $X\times\{t\}$. Indeed, the
restriction is $S^0+\llb 0,t\rrb\times T$.

The slicing formula
then gives
\[
  \partial\langle W,\vartheta,t\rangle
  =-\partial(U\llcorner\{\vartheta<t\})=-T^t.
\]
The slice is supported in $X\times\{t\}=\vartheta^{-1}(t)$.
By the coarea inequality, some such slice has mass at most $2\M(W)$.
Identification 
with $X$ gives
\[
  \operatorname{Fillvol}_X(T)
  \leq C_d\bigl(s+\M(T)\bigr)^{(d+2)/(d+1)}.
\]

We apply this same unit-width estimate after multiplying the
metrics of $X$ and $E$ by $\lambda>0$. Rescaling the masses
back gives
\[
  \operatorname{Fillvol}_X(T)
  \leq C_d\lambda^{-1/(d+1)}
        \bigl(\lambda s+\M(T)\bigr)^{(d+2)/(d+1)}.
\]
For $T\neq0$ we have $s>0$. Set $\lambda=\M(T)/s$ and take
the infimum over the ambient fillings $S$. The case $T=0$
is immediate.
\end{proof}

The next lemma extends a prescribed Hausdorff convergence to
realize convergence of representatives of points in a separable
subset of the ultralimit.

\begin{lemma}\label{lem:relative-ultralimit-realization}
Let $(X_i,x_i)$ be pointed metric spaces and let $K_i\subset X_i$
be compact subsets containing $x_i$. Suppose isometric embeddings
$\varphi_i:K_i\to Z$ into a metric space $Z$ satisfy
$\varphi_i(K_i)\to K$ in Hausdorff distance, where $K\subset Z$
is compact. Let $X_\omega$ be the pointed ultralimit of the $X_i$
with respect to a nonprincipal ultrafilter $\omega$ and $\iota_\omega:K\to X_\omega$ 
the induced isometric embedding.

For every separable subset $A\subset X_\omega$, there is a metric
space $W$ containing isometric copies of $Z$, $X_\omega$, and all
$X_i$, with the copies of each $K_i$ identified via $\varphi_i$
and the copies of $K$ identified via $\iota_\omega$. For each
$p\in\overline A$ and every representative $(p_i)$ of $p$, one has
$\lim_\omega d_W(p_i,p)= 0$.
\end{lemma}

\begin{proof}
For $x,y\in K$, choose $a_i,b_i\in K_i$ with
$\varphi_i(a_i)\to x$ and $\varphi_i(b_i)\to y$. Then
$$
 d_{X_\omega}((a_i),(b_i))
 =\lim_\omega d_Z(\varphi_i(a_i),\varphi_i(b_i))=d_Z(x,y).
$$
This shows that $\iota_\omega$ is well defined and isometric.
For any sequence $a_i\in K_i$, Hausdorff convergence gives
$d_Z(\varphi_i(a_i),K)\to0$. By compactness of $K$, the sequence
$\varphi_i(a_i)$ has an $\omega$-limit in $K$, and $(a_i)$
is the image of this limit under $\iota_\omega$
(cf.~\cite[Lemma~6.5]{WengerAlmostEuclidean}).

Glue $Z$ and $X_\omega$ along $K$ via $\iota_\omega$ and denote
the resulting space by $W_0$. Attaching $X_i$ to $W_0$ via
$\varphi_i$ identifies each $a\in K_i$ with its prescribed image
$\varphi_i(a)\in Z$. With these identifications alone, the
distance from $p\in X_\omega$ to $X_i$ is at least
$d_{X_\omega}(p,\iota_\omega(K))$, since $X_\omega$ is attached
to $Z$ only along $K$. To realize convergence of representatives
for points of $\overline A$, we add edges between selected points
of $X_\omega$ and their representatives in $X_i$.

Choose a sequence $(p^r)_{r\ge1}$ in $X_\omega$ whose closure
contains $A$, and for each $r$ choose a representative
$(p_i^r)_i$ of $p^r$. To preserve the metrics on $X_i$ and $W_0$
when adding edges, we use the following estimate for each fixed $r$:
$$
 e_i^r:=\sup_{a\in K_i}
 \left|d_{X_i}(p_i^r,a)-d_{W_0}(p^r,\varphi_i(a))\right|
 \longrightarrow_\omega 0.
$$
Choose a finite $\varepsilon$-net in $K$. For each of its points
$k$, choose $k_i\in K_i$ with $\varphi_i(k_i)\to k$. At $a=k_i$,
both distances  converge along $\omega$ to
$d_{X_\omega}(p^r,\iota_\omega(k))$. For large $i$ these points
$k_i$ form a $3\varepsilon$-net in $K_i$. Since both distance
functions are $1$-Lipschitz on $K_i$, this gives
$\limsup_\omega e_i^r\le 6\varepsilon$. Letting
$\varepsilon\to0$ proves the estimate.

For each $i$ and $r$, set
$$
 \ell_i^r=\frac1i+\max\left\{e_i^r,\
 \max_{1\le s\le r}
 \left|d_{X_i}(p_i^r,p_i^s)-d_{X_\omega}(p^r,p^s)\right|\right\}.
$$
Then $\lim_\omega\ell_i^r=0$ for every fixed $r$.

For each $i$, add an edge of length $\ell_i^r$ from $p_i^r$ to
$p^r$ for every $r$. Define a pseudo-metric on $X_i\sqcup W_0$
by taking the infimum of the total cost of finite sequences of
steps. A step within either space costs the distance in that
space, crossing an edge costs its length, and passing between
$a\in K_i$ and $\varphi_i(a)\in W_0$ costs zero. This definition
of the cross-distance is analogous to the one used in the proof
of \cite[Theorem~1.6]{IkonenWenger}.

For two edges, the discrepancy between the distances of their
endpoints in the two spaces is bounded by the length of the
edge with larger index. For an edge from $p_i^r$ to $p^r$ and
an identification along $K_i$, the discrepancy is at most $e_i^r$.
For two identifications it is zero. In all cases the discrepancy
is at most the sum of the two crossing costs. Thus an excursion
through the other space costs at least the distance between
its endpoints in the original space. Consequently $X_i$ and
$W_0$ embed isometrically in the metric quotient. Glue these
quotients along their common copy of $W_0$ to obtain $W$.

For $p=(p_i)\in\overline A$ and each $r$, we have
$$
 d_W(p_i,p)\le d_{X_i}(p_i,p_i^r)+\ell_i^r+d_{X_\omega}(p^r,p).
$$
The limit superior along $\omega$ is at most
$2d_{X_\omega}(p^r,p)$, and density gives $d_W(p_i,p)\to_\omega0$.
In particular, $\lim_\omega d_W(p,X_i)=0$ for every
$p\in\overline A$. 
% The functions $p\mapsto d_W(p,X_i)$ are
% $1$-Lipschitz, so a finite net in $C$ gives uniform convergence
% on $C$ along $\omega$.
\end{proof}

\begin{proof}[Proof of Theorem~\ref{thm:ultralimit-lifting}]
By Theorem~\ref{thm:main}, approximate $S$ in normal mass by
Lipschitz chains $\alpha_j$ with $\Vol(\alpha_j)\to\M(S)$ and
$\Vol(\partial\alpha_j)\to\M(\partial S)$. The union $A$ of their
carriers is separable, since each chain consists of finitely many
simplices with compact images. Apply
Lemma~\ref{lem:relative-ultralimit-realization} to $A$ and embed
$W$ isometrically into $E=\ell^\infty(W)$ by the Kuratowski
embedding. We have $\lim_\omega d_E(p,X_i)=0$ uniformly on the
carrier of each $\alpha_j$. Since the inclusions agree on $K_i$
and $K$, pushing forward from $Z$ to $E$ gives
$\lim_\omega\mathcal F_E(B_i-\partial S)=0$. We identify the
currents with their images in $E$ in the remainder of the proof.

For each fixed $j$, lift the simplices in order of increasing
dimension by Lemma~\ref{lem:lift-simplex}. Use the same lift for each
occurrence of a given simplex. The face lifts agree on intersections
and give uniformly Lipschitz boundary data. For edges this follows from convergence of the lifted
endpoints. 
Let $\alpha_{j,i}$ be the lifted chain.
By Lemma~\ref{lem:lift-simplex},
$\limsup_\omega\Vol(\alpha_{j,i})\leq\Vol(\alpha_j)$ and
$\limsup_\omega\Vol(\partial\alpha_{j,i})\leq\Vol(\partial\alpha_j)$.
The straight-line homotopy in $E$ gives
$\lim_\omega\mathcal F_E([\alpha_{j,i}]-[\alpha_j])=0$
by \eqref{eq:prism}. A diagonal argument gives integral currents $P_i$ in
$X_i$ with
$$
  \lim_\omega\mathcal F_E(P_i-S)=0,\qquad
  \limsup_\omega\M(P_i)\leq\M(S),\qquad
  \sup_i\M(\partial P_i)<\infty.
$$

Set $D_i=B_i-\partial P_i$. These cycles have uniformly bounded
mass and $\lim_\omega \mathcal F_E(D_i)=0$, since the boundary operator
is $1$-Lipschitz for the flat norm. Write $D_i=U_i+\partial V_i$
in $E$ with $\lim_\omega(\M(U_i)+\M(V_i))=0$. If $m>1$, the Euclidean
isoperimetric inequality in $E$ fills the cycles $U_i$ with mass
tending to zero along $\omega$. If $m=1$, then $U_i=0$ whenever
$\M(U_i)<1$. Thus $\lim_\omega\operatorname{Fillvol}_E(D_i)=0$.
Lemma~\ref{lem:intrinsic-small-fillings} now gives
$\lim_\omega\operatorname{Fillvol}_{X_i}(D_i)=0$.

Each $D_i$ has a filling in $X_i$ by the isoperimetric inequality. 
Choose such
fillings $R_i$ with $\lim_\omega \M(R_i)=0$.
Then $S_i=P_i+R_i$ has boundary $B_i$ and satisfies
$\lim_\omega\mathcal F_E(S_i-S)=0$. The upper mass bound for $P_i$
and lower semicontinuity of mass give $\lim_\omega\M(S_i)=\M(S)$.
\end{proof}

\begin{proof}[Proof of Corollary~\ref{cor:ultralimit-minimizers}]
Fix a nonprincipal ultrafilter $\omega$. For any $S\in\I_m(X_\omega)$
with $\partial S=(\iota_\omega)_\#\partial T$, apply
Theorem~\ref{thm:ultralimit-lifting} with $B_i=\partial T_i$ and
$B=\partial T$. Lower semicontinuity in $Z$ and minimality in $X_i$
give
$$
  \M(T)\leq\liminf_\omega\M(T_i)
  \leq\limsup_\omega\M(T_i)\leq\M(S).
$$
Taking $S=(\iota_\omega)_\#T$ gives $\lim_\omega\M(T_i)=\M(T)$.
Since this holds for every nonprincipal $\omega$, the masses
converge to $\M(T)$ in the ordinary sense.
\end{proof}

% \begin{example}\label{ex:property-L}
% There is a compact CAT(0) space without local property~L.
% \end{example}

% \begin{proof}
% By \cite[Theorem~4]{NaorRabani}, take finite metric spaces $Y_m$, subsets
% $A_m\subset Y_m$, and $1$-Lipschitz maps $f_m:A_m\to H_m$ into
% finite-dimensional Hilbert spaces such that every extension to $Y_m$ has
% Lipschitz constant at least $m$. Rescale both metrics so that
% $\operatorname{diam}Y_m\leq2^{-m}$, and translate $f_m$ to send a chosen
% $a_m\in A_m$ to $0$. Identify the $a_m$ to a point $o$ in a metric
% wedge $Y$, with distances between different summands given by the sum of
% their distances to $o$. Then $A=\bigcup_m A_m$ is compact.

% In the Hilbert direct sum of the $H_m$, let $K$ be the closed convex hull
% of the union of the images $f_m(A_m)$. This union is compact, since the
% images tend uniformly to $0$. Its closed convex hull is compact as well,
% by taking convex hulls of finite nets. With its Hilbert metric $K$ is
% CAT(0). The maps $f_m$ give a $1$-Lipschitz map $f:A\to K$, since
% $\|f_m(x)\|\leq d(x,o)$ and distances across summands add. Every
% $\varepsilon$-neighborhood of $A$ in $Y$ contains all $Y_m$ with
% $2^{-m}<\varepsilon$. Composing an extension with the orthogonal
% projections onto $H_m$ would give extensions of these $f_m$ with one
% Lipschitz bound. The same argument on a tail applies to any neighborhood of
% $0$ in $K$, so $K$ fails local property~L at $0$.
% \end{proof}

\begingroup
\raggedright

\endgroup
\end{document}